\documentclass{amsart}[12pt]
\usepackage{a4}
\usepackage[T1]{fontenc}
\usepackage[latin1]{inputenc}
\usepackage{hyperref}
\usepackage[dvips]{graphics}
\usepackage{amsmath}
\usepackage{amssymb}
\usepackage{amsopn}
\usepackage{amsbsy}
\usepackage{amstext}
\usepackage{latexsym}
\usepackage{amsfonts}
\usepackage{array}
\usepackage{epsfig}
\usepackage{mathrsfs}
\title[Formulas giving prime numbers under Cramér's conjecture]{Formulas giving prime numbers under Cramér's conjecture}
\author[B. FARHI]{Bakir FARHI}
\date{\bf November 24\textsuperscript{th}, 2006}

\newtheorem{thm}{Theorem}

\newtheorem{coll}[thm]{Corollary}
\newtheorem{conj}[thm]{Conjecture}

\let\epsilon=\varepsilon

\def\EMdash{\leavevmode\hbox to 7.5mm{\vrule height .63ex depth -.59ex
    width 5.4mm\hfill}}

\begin{document}
\maketitle \baselineskip=6mm
\begin{center}
Département de Mathématiques, Université du Maine, \\
Avenue Olivier Messiaen, 72085 Le Mans Cedex 9, France. \\
Bakir.Farhi@univ-lemans.fr
\end{center}
\begin{abstract}
Under Cramer's conjecture concerning the prime numbers, we prove
that for any $\xi > 1$, there exists a real $A = A(\xi)
> 1$ for which the formula $[A^{n^{\xi}}]$ (where $[.]$ denotes
the integer part) gives a prime number for any positive integer
$n$. Under the same conjecture, we also prove that for any
$\varepsilon > 0$, there exists a real $B = B(\varepsilon) > 0$
for which the formula $[B . {n!}^{2 + \epsilon}]$ gives a prime
number for any sufficiently large positive integer $n$.
\end{abstract}
\hfill{--------} 
\section{Introduction}
Throughout this article, we will let $[x]$ denote the integer part
of a giving real number $x$; also, we will let ${(p_n)}_{n \in
\mathbb{N}}$ denote the sequence of all prime numbers and we will
put $\Delta p_n := p_{n + 1} - p_n$ for any $n \in \mathbb{N}$.
Further, if $A$ is a subset of $\mathbb{R}$ and $x$ is a real
number, we will let $A + x$ denote the subset of $\mathbb{R}$
defined by: $A + x := \{a + x | a \in A\}$.

In \cite{mi}, W. H. Mills proved the existence of an absolute
constant $A > 1$ for which $[A^{3^n}]$ is a prime number for any
positive integer $n$ and in \cite{w}, E. M. Wright proved the
existence of an absolute constant $\alpha > 0$ for which the
infinite sequence $[\alpha] , [2^{\alpha}] , [2^{2^{\alpha}}] ,
\dots $ is composed of prime numbers. Let us describe the method
used by these two authors. They start from an upper bound of
$\Delta p_n$ as a function of $p_n$. Such upper bound allows to
construct an increasing function $h$ (more or less elementary,
according to the used upper bound of $\Delta p_n$) such that
between any two consecutive terms of the sequence ${(h(n))}_n$,
there is at least one prime number. Setting $f_n := h \circ \dots
\circ h$ (where $h$ is applied $n$ times), they deduce from the
last fact, the existence of a real constant $A$ for which the
sequence ${([f_n(A)])}_n$ is consisted of prime numbers.

By this method, Wright uses the upper bound $\Delta p_n \leq p_n$
which is nothing else than Bertrand's postulate and Mills uses
Ingham's upper bound: $\Delta p_n \leq p_{n}^{\frac{5}{8} +
\varepsilon}$ (which holds for any sufficiently large $n$ compared
to a given $\epsilon > 0$). The functions $h$ which derive from
these upper bounds are $h(x) = 2^x$ for Wright and $h(x) = x^3$
for Mills. Then the Theorems of \cite{mi} and \cite{w} follow.

Notice that the more the upper bound of $\Delta p_n$ is refined,
the more the function $h$ is small and the more the obtained
sequence of prime numbers grows slowly (remark for instance that
the sequence of Mills grows slower than the Wright's one). From
this fact, in order to have a sequence of prime numbers which
grows more slow again, we must use upper bounds more refined for
$\Delta p_n$. But up to now even the powerful Riemann hypothesis
gives only the estimate $\Delta p_n = O(p_n^{1/2} \log p_n)$. A
famous conjecture (which is little too strong than this last
estimate) states that between two consecutive squares, there is
always a prime number (see \cite{guy}). So, according to this
conjecture, the function $h(x) = x^2$ is admissible in the method
described above, which we permits to conclude the existence of a
constant $B > 1$ for which $[B^{2^n}]$ is a prime number for any
positive integer $n$. We thus obtain (via this conjecture), a
sequence of prime numbers growing slower than the Mills'one.

By leaning on heuristic and probabilistic arguments, H. Cramér
\cite{cra} was leaded to conjecture that we have $\Delta p_n =
O(\log^2 p_n)$; further it's known that $\Delta p_n =
O(\log{p_n})$ cannot hold (see \cite{we}). Thus, by taking in the
method described above $h(x) = c \log^2 x$ $(c > 0)$, we obtain
(via Cramér's conjecture) sequences of prime numbers having
explicit form and growing much slower than the Mills'one. The
inconvenient in this application is that the explicit form in
question $[f_n(A)]$ is not elementary, because $f_n$ doesn't have
(in this case) a simple expression as function of $n$.

To cope with this problem, we were leaded to generalize
Mills'method by considering instead of one function $h$, a
sequence of functions ${(h_m)}_m$ and in this situation $f_n$ is
rather the composition of $n$ functions $h_0 , \dots , h_{n - 1}$.
This allows fundamentally to give for $f_n$ the form which we
want, then if we set $h_n := f_{n + 1} \circ f_{n}^{- 1}$, we have
only to check whether it's true that for any $n$ and any $x$
sufficiently large (relative to $n$), the interval $[h_n(x) ,
h_n(x + 1) - 1[$ contains at least one prime number or not. In the
affirmative case, we will deduce the existence of a real $A$ for
which the formula $[f_n(A)]$ gives a prime number for any positive
integer $n$ (see Theorem \ref{t1} and its proof).

Under a conjecture less strong than the Cramér's one, we derive
from this generalization two new types of explicit formulas giving
prime numbers. We also give other applications of our main result
(outside the subject of prime numbers) and we conclude this
article by some open questions related to the results which we
obtain.
\section{Results}
The main result of this article is the following
\begin{thm}\label{t1}
Let $I = ]a , b[$ (with $a , b \in \overline{\mathbb{R}}$, $a <
b$) be an open interval of $\mathbb{R}$, $n_0$ be a non-negative
integer and ${(f_n)}_{n \geq n_0}$ be a sequence of real functions which are differentiable and increasing on $I$.\\
We assume that the functions $\frac{f'_{n + 1}}{f'_n}$ $(n \geq
n_0)$ are nondecreasing on $I$ and that for all $x \in I$, the
numerical sequence ${(f_n(x))}_{n \geq n_0}$ is increasing. \\
We also assume that there exists a real function $g$,
nondecreasing on $\mathbb{R}$ and verifying:
\begin{equation}\label{eq1}
g \circ f_{n + 1} (x) ~\leq~ \frac{f'_{n + 1}}{f'_n}(x)
\hspace{2cm} \text{\rm($\forall n \geq n_0$ and $\forall x \in
I$).}
\end{equation}
Then, for any sequence of integers ${(u_n)}_{n \in \mathbb{N}}$,
verifying: $\displaystyle \limsup_{n \rightarrow + \infty} u_n = +
\infty$,
\begin{equation}\label{eq2}
u_{n + 1} - u_n ~\leq~ g(u_n) - 1 \hspace{2cm} \text{\rm(for all
$n$ sufficiently large)}
\end{equation}
and whose one at least of the terms $u_n$ ($n$ satisfying {\rm
(\ref{eq2})}) belongs to $f_{n_0}(I) \cap (f_{n_0}(I) - 1)$, there
exists a real $A \in I$ for which the sequence ${([f_n(A)])}_{n
\geq n_0}$ is an increasing subsequence of ${(u_n)}_n$.
\end{thm}
\vskip 1mm \noindent {\bf Proof:} By shifting if necessary the
sequence of functions ${(f_n)}_{n \geq n_0}$, we may assume that
$n_0 = 0$ and by shifting if necessary the sequence ${(u_n)}_n$,
we may assume that we have (more generally than (\ref{eq2})):
\begin{equation}\label{eq2'}
u_{n + 1} - u_n ~\leq~ g(u_n) - 1 \hspace{2cm} \text{(for all $n
\in \mathbb{N}$).} \tag{$2'$}
\end{equation}
We begin this proof by some remarks and preliminary notations
which allow to lighten better the situation of the Theorem.\\
Giving $n \in \mathbb{N}$, since the function $f_n$ is assumed
differentiable (so continuous) and increasing on $I = ]a , b[$,
then it's a bijection from $I$ into $\displaystyle f_n(I) =
]\lim_{x \rightarrow a} f_n(x) , \lim_{x \rightarrow b} f_n(x)[ =
]\lambda_n , \mu_n[$, where $\displaystyle \lambda_n := \lim_{x
\rightarrow a} f_n(x)$ and $\displaystyle \mu_n := \lim_{x
\rightarrow b} f_n(x)$
($\lambda_n$ and $\mu_n$ belong to $\overline{\mathbb{R}}$). \\
Now, let us introduce the following functions:
\begin{equation*}
\begin{split}
h_n :~ ]\lambda_n , \mu_n[ ~\longrightarrow~ ]\lambda_{n + 1} , \mu_{n + 1}[ \\
h_n := f_{n + 1} \circ f_{n}^{-1}
\end{split} ~~~~~~(\forall n \in \mathbb{N}).
\end{equation*}
Giving $n \in \mathbb{N}$, since (from the hypothesis of the
Theorem), the functions $f_n$ and $f_{n + 1}$ are differentiable
and increasing on $I$, then the function $h_n$ is differentiable
and increasing on $]\lambda_n , \mu_n[$. Further, the hypothesis
of the Theorem concerning the increase of the numerical sequences
${(f_n(x))}_n$ $(x \in I)$ amounts to:
\begin{equation}\label{eq3}
h_n(x) ~>~ x ~~~~~~ \text{($\forall n \in \mathbb{N}$ and $\forall
x \in ]\lambda_n , \mu_n[$).}
\end{equation}
Next, let us show that for any $n \in \mathbb{N}$, the function
$h_n$ is convex on $]\lambda_n , \mu_n[$. To do this, we check
that the derivative $h'_n$ of each function $h_n$ $(n \in
\mathbb{N})$ is nondecreasing on the interval $]\lambda_n ,
\mu_n[$. Giving $n \in \mathbb{N}$, we have:
$$h'_n = (f_{n + 1} \circ f_{n}^{-1})' = (f_{n}^{-1})' . f'_{n + 1} \circ f_{n}^{-1} = \frac{f'_{n + 1}
\circ f_{n}^{-1}}{f'_n \circ f_{n}^{-1}} = \left(\frac{f'_{n +
1}}{f'_n}\right) \circ f_{n}^{-1} .$$ Since (from the hypothesis
of the Theorem), the function $\frac{f'_{n + 1}}{f'_n}$ is
nondecreasing on $I$ and the function $f_{n}^{-1}$ is increasing
on $f_n(I) = ]\lambda_n , \mu_n[$ (because $f_n$ is increasing on
$I$), then (as a composite of two nondecreasing functions), the
function $h'_n$ is nondecreasing on $]\lambda_n , \mu_n[$. So the
function $h_n$ is effectively convex on
$]\lambda_n , \mu_n[$.\\
The rest of this proof consists of the three following steps:\\
\underline{1\textsuperscript{st} Step:}~\vspace{1mm} We are going
to show that we have:
\begin{equation}\label{eq4}
g \circ h_n(y) ~\leq~ h_n(y + 1) - h_n(y) ~~~~~~ \text{($\forall n
\in \mathbb{N}$ and $\forall y \in ]\lambda_n , \mu_n - 1[$).}
\end{equation}
(We will see further that the interval $]\lambda_n , \mu_n - 1[$
is never empty).\\
Let $n \in \mathbb{N}$ and $y \in ]\lambda_n , \mu_n - 1[$ fixed and set $x := f_{n}^{-1}(y)$.\\
The convexity of $h_n$ on $]\lambda_n , \mu_n[$, proved above,
implies that we have:
$$h_n(u) ~\geq~ h'_n(t) (u - t) + h_n(t) \hspace{2cm} \text{(for all $t , u \in ]\lambda_n , \mu_n[$).}$$
By taking in this last inequality $t = y$ and $u = y + 1$, we
obtain:
\begin{eqnarray*}
h_n(y + 1) - h_n(y) & \geq & h'_n(y) \\
& = & \left(\frac{f'_{n + 1}}{f'_n}\right)(x) ~~ \text{(because
$h'_n = \frac{f'_{n + 1}}{f'_n} \circ f_{n}^{-1}$
and $x = f_{n}^{-1}(y)$)} \\
& \geq & g \circ f_{n + 1} (x) ~~~~ \text{(from the hypothesis
(\ref{eq1}) of the Theorem)} \\
& = & g \circ f_{n + 1} \circ f_{n}^{-1}(y) \\
& = & g \circ h_n(y) .
\end{eqnarray*}
The relation (\ref{eq4}) follows.\\
\underline{2\textsuperscript{nd} Step:} We are going to construct
an increasing sequence ${(k_n)}_{n \in \mathbb{N}}$ of
non-negative integers such that the subsequence of ${(u_n)}_n$
with general term $v_n = u_{k_n}$ satisfies:
\begin{equation*}
\left\{
\begin{split}
&v_n \in ]\lambda_n , \mu_n - 1[ \\
&\text{and}\\ &h_n(v_n) ~\leq~ v_{n + 1} ~<~ h_n(v_n + 1) - 1
\end{split} ~~~~~~ (\forall n \in \mathbb{N}).\right.
\end{equation*}
We proceed by induction as follows:\\
$\bullet$ We pick $k_0 \in \mathbb{N}$ such that $u_{k_0} \in
f_0(I) \cap (f_0(I) - 1) = ]\lambda_0 , \mu_0 - 1[$. Notice that
the existence of such integer $k_0$ is an hypothesis of the Theorem.\\
$\bullet$ If, for some $n \in \mathbb{N}$, an integer $k_n \in
\mathbb{N}$ is chosen such that $u_{k_n} \in ]\lambda_n , \mu_n -
1[$, let:
$$X_n ~:=~ \left\{k \in \mathbb{N} ~|~ k > k_n ~\text{and}~ u_k \geq h_n(u_{k_n})\right\} .$$
From the hypothesis $\limsup_{n \rightarrow + \infty} u_n = +
\infty$, this subset $X_n$ of $\mathbb{N}$ is nonempty, it thus
admits a smallest element which we take as the choice of $k_{n +
1}$. So, we have:
$$k_{n + 1} > k_n ~,~ u_{k_{n + 1}} ~\geq~ h_n(u_{k_n}) ~~\text{and}~~ k_{n + 1} - 1 \not\in X_n .$$
We claim that the facts ``$k_{n + 1} > k_n$'' and ``$k_{n + 1} - 1
\not\in X_n$'' imply:
\begin{equation}\label{eq9}
u_{k_{n + 1} - 1} ~<~ h_n(u_{k_n}) .
\end{equation}
Indeed: either $k_{n + 1} = k_n + 1$ in which case we have
$u_{k_{n + 1} - 1} = u_{k_n} < h_n(u_{k_n})$ (from (\ref{eq3}));
or $k_{n + 1} > k_n + 1$, that is $k_{n + 1} - 1 > k_n$, but since
$k_{n + 1} - 1 \not\in X_n$, we are forced to have (also in this case) $u_{k_{n + 1} - 1} < h_n(u_{k_n})$ as required.\\
It follows that we have:
\begin{eqnarray*}
u_{k_{n + 1}} & \leq & u_{k_{n + 1} - 1} + g(u_{k_{n + 1} - 1}) -
1 \hspace{6mm} \text{(from (\ref{eq2'}))} \\
& < & h_n(u_{k_n}) + g \circ h_n(u_{k_n}) - 1 \hspace{5mm} \text{(by using (\ref{eq9}) and the non-decrease of $g$)} \\
& \leq & h_n(u_{k_n} + 1) - 1 \hspace{2cm} \text{(from
(\ref{eq4}))} .
\end{eqnarray*}
Hence:
$$u_{k_{n + 1}} ~<~ h_n(u_{k_n} + 1) - 1 .$$
We thus have:
$$h_n(u_{k_n}) ~\leq~ u_{k_{n + 1}} ~<~ h_n(u_{k_n} + 1) - 1 .$$
Finally, since the function $h_n$ takes its values in $]\lambda_{n
+ 1} , \mu_{n + 1}[$, then the last double inequality does show
that $u_{k_{n + 1}} \in ]\lambda_{n + 1} , \mu_{n + 1} - 1[$. This
ensures the good working of this induction process which gives the required sequence ${(k_n)}_n$.\\
Notice also that the subsequence ${(v_n)}_n$ of ${(u_n)}_n$ which
we have just constructed is increasing, because for any $n \in
\mathbb{N}$, we have: $v_{n + 1} \geq h_n(v_n)
> v_n$ (from (\ref{eq3})).\\
\underline{3\textsuperscript{rd} step (conclusion):} To conclude
this proof, we will show the existence of a real $A \in I$ for
which we have for any $n \in \mathbb{N}$:
$v_n = [f_n(A)]$.\\
To do this, we introduce two real sequences ${(x_n)}_n$ and
${(y_n)}_n$, with terms in $I$, which we define by:
$$x_n := f_{n}^{-1}(v_n) ~~\text{and}~~ y_n := f_{n}^{-1}(v_n + 1) ~~~~~~ (\forall n \in \mathbb{N}) .$$
Since the functions $f_n$ are increasing, we have $x_n < y_n$ for
all $n \in \mathbb{N}$. We claim that the sequence ${(x_n)}_n$ is
nondecreasing and that the sequence ${(y_n)}_n$ is decreasing.
Indeed, for any $n \in \mathbb{N}$, we have:
$$x_n = f_{n}^{-1}(v_n) = f_{n + 1}^{-1} \circ h_n(v_n) \leq f_{n + 1}^{-1}(v_{n + 1}) = x_{n + 1}$$
and
$$y_n = f_{n}^{-1}(v_n + 1) = f_{n + 1}^{-1} \circ h_n(v_n + 1) > f_{n + 1}^{-1}(v_{n + 1} + 1) = y_{n + 1} .$$
(Where in these last relations, we have just use the facts that
$f_{n + 1}^{-1}$ is increasing and $h_n(v_n) \leq
v_{n + 1} < h_n(v_n + 1) - 1$).\\
The $[x_n , y_n]$ $(n \in \mathbb{N})$ are thus nested inclosed
intervals of $\mathbb{R}$. Consequently their intersection is
nonempty (according to the Cantor's intersection theorem). Pick
$A$ an arbitrary real belonging to this intersection, that is:
$x_n \leq A \leq y_n$ for all $n \in \mathbb{N}$ (in particular $A
\in I$). In fact $A$ verifies even:
$$x_n ~\leq~ A ~<~ y_n ~~~~~~ (\forall n \in \mathbb{N}) ,$$
because if $A = y_m$ for some $m \in \mathbb{N}$, we will have
(from the decreasing of the sequence ${(y_n)}_n$):
$A > y_{m + 1}$, contradicting the inequality $A \leq y_{m + 1}$.\\
It follows from the increase of the functions $f_n$ that we have:
$$f_n(x_n) ~\leq~ f_n(A) ~<~ f_n(y_n) ~~~~~~ (\forall n \in \mathbb{N}) ,$$
that is:
$$v_n ~\leq~ f_n(A) ~<~ v_n + 1 ~~~~~~ (\forall n \in \mathbb{N}) .$$
Then (since $v_n$ is an integer for all $n \in \mathbb{N}$):
$$[f_n(A)] ~=~ v_n ~~~~~~ (\forall n \in \mathbb{N}) .$$
This completes the
proof.\penalty-20\null\hfill$\blacksquare$\par\medbreak \noindent
{\bf Remarks:} Mills'theorem \cite{mi} can be find again by
applying the above Theorem \ref{t1} for $I = ]1 , + \infty[$, $n_0
= 0$, $f_n(x) = x^{3^n}$ ($\forall n \in \mathbb{N}$ and $\forall
x \in I$), $g(x) = x^{2/3}$ if $x > 0$ and $g(x) = 0$ if $x \leq
0$ and ${(u_n)}_n$ the sequence of the prime numbers. In this
application, we check the relation (\ref{eq1}) of Theorem \ref{t1}
by a simple calculus and we deduce the relation (\ref{eq2}) of the
same Theorem from the Ingham's estimate quoted in the
introduction. The remaining hypothesis of Theorem \ref{t1} are
immediately verified.

The Wright's theorem \cite{w} can be also find again by applying
Theorem \ref{t1} for $I = ]0 , + \infty[$, $n_0 = 0$, ${(f_n)}_{n
\in \mathbb{N}}$ the sequence of functions which is defined on $I$
by: $f_0 = {\rm{Id}}_I$ and $f_{n + 1} = 2^{f_n}$ $(\forall n \in
\mathbb{N})$, $g(x) = (\log{2}) x$ $(\forall x \in \mathbb{R})$
and ${(u_n)}_n$ the sequence of the prime numbers. In order to
check the relation (\ref{eq1}) of Theorem \ref{t1}, remark that we
have for any $n \in \mathbb{N}$: $\frac{f'_{n + 1}}{f'_n} =
(\log{2}) f_{n + 1}$. As for the relation (\ref{eq2}), it's (in
this application) a consequence of the prime numbers theorem; but
it can be obtained by using Chebyshev's elementary arguments (see
\cite{har}). The remaining hypothesis of Theorem
\ref{t1} are immediately verified.~\vspace{1mm}\\
{\bf{N. B.---}} It's interesting to remark that in the two above
applications of Theorem \ref{t1}, the sequence of functions
${(h_n)}_n$ introduced in the proof of the latter $(h_n := f_{n +
1} \circ f_{n}^{-1})$ is constant. Indeed, for the first
application, we have $h_n(x) = x^3$ $(\forall n \in \mathbb{N})$
and for the second one, we have $h_n(x) = 2^x$ $(\forall n \in
\mathbb{N})$. As explained in the introduction, the fact to be
able to take ${(h_n)}_n$ not constant is the crucial point of our
approach. In the following, we are going to give some applications
of Theorem \ref{t1} in which this sequence ${(h_n)}_n$ is not
constant. If we admit the following Conjecture (which is less
strong than the Cramér's one \cite{cra}), we obtain two new types
of explicit sequences of prime numbers, which grow much more slow
than the ones of Mills and Wright.
\begin{conj}\label{conj1}
there exists an absolute constant $k > 1$ such that:
$$\Delta p_n ~=~ O\left((\log{p_n})^k\right) .$$
\end{conj}
Under this Conjecture, we obtain by applying Theorem \ref{t1}, the
two following Corollaries:
\begin{coll}\label{coll1}
Under the above Conjecture \ref{conj1}, there exists for all real
$\xi > 1$, a real $A = A(\xi) > 1$ for which the sequence
${([A^{n^{\xi}}])}_{n \geq 1}$ is an increasing sequence of prime
numbers.
\end{coll}
\vskip 1mm \noindent {\bf Proof:} Let $\xi > 1$ be fixed, $k
> 1$ be an admissible constant in Conjecture \ref{conj1} and $a
> 1$ be a real such that:
\begin{eqnarray}
\left(\log{x}\right)^{k + 1} & \leq & x^{1/2} \hspace{2cm} (\forall x > a) \label{eq5} \\
(n + 1)^{k + 1} & \leq & a^{\frac{1}{2} n^{\xi - 1}}
\hspace{1.6cm} (\forall n \geq 1). \label{eq6}
\end{eqnarray}
($a$ exists because $\lim_{x \rightarrow + \infty}
\frac{(\log{x})^{k + 1}}{x^{1/2}} = 0$ and $\lim_{n \rightarrow + \infty} (n + 1)^{\frac{2 (k + 1)}{n^{\xi - 1}}} = 1$).\\
We apply Theorem \ref{t1} for $I = ]a , + \infty[$, $n_0 = 1$,
$f_n(x) = x^{n^{\xi}}$ $(\forall n \geq 1, \forall x \in I)$,
$g(x) = (\log{x})^{k + 1}$ if $x > 1$ and $g(x) = 0$ if $x \leq 1$
and ${(u_n)}_n$ the sequence of the prime numbers.
Let us check the hypothesis of Theorem \ref{t1}:\\
The functions $f_n$ are clearly increasing and differentiable on
$I$. We have $f'_n(x) = n^{\xi} x^{n^{\xi} - 1}$ ($\forall n \geq
1$, $\forall x \in I$), then: $\frac{f'_{n + 1}}{f'_n}(x) =
(\frac{n + 1}{n})^{\xi} x^{(n + 1)^{\xi} - n^{\xi}}$ ($\forall n
\geq 1$, $\forall x \in I$). We thus see that the functions
$\frac{f'_{n + 1}}{f'_n}$ $(n \geq 1)$ are nondecreasing on $I$.
Further, if $x$ is a fixed real in $I$, the numerical sequence
${(f_n(x))}_{n \geq 1}$ is clearly increasing. Now, we have for
any integer $n \geq 1$ and for any real $x \in I$:
\begin{eqnarray*}
g \circ f_{n + 1}(x) & = & (n + 1)^{\xi (k + 1)} (\log{x})^{k + 1}
\\
& \leq & a^{\frac{1}{2} \xi n^{\xi - 1}} x^{1/2} \hspace{2cm}
\text{(from (\ref{eq5}) and (\ref{eq6}))} \\
& \leq & x^{\xi n^{\xi - 1}} \hspace{2.8cm} \text{(because
$x > a$ and $\xi n^{\xi - 1} > 1$)} \\
& \leq & x^{(n + 1)^{\xi} - n^{\xi}} \hspace{2.2cm} \text{(because
$\xi n^{\xi - 1} \leq (n +
1)^{\xi} - n^{\xi}$)} \\
& \leq & \frac{f'_{n + 1}}{f'_n}(x) .
\end{eqnarray*}
The relation (\ref{eq1}) of Theorem \ref{t1} follows. \\
Next, the relation (\ref{eq2}) of Theorem \ref{t1} (related to
this application) follows immediately from the conjecture
\ref{conj1} (admitted in this context). Finally, $f_{n_0}(I) \cap
(f_{n_0}(I) - 1) = ]a , + \infty[$ does contains prime numbers as
large as we want. The hypothesis of Theorem \ref{t1} are thus all
satisfied, so we can apply this latter to the present situation.
The corollary \ref{coll1} follows from this
application.\penalty-20\null\hfill$\blacksquare$\par\medbreak
\begin{coll}\label{coll2}
Assume that Conjecture \ref{conj1} is true and let $k > 1$ be an
admissible constant in this Conjecture. Then, for any positive
real $\varepsilon$, there exists an integer $n_0 = n_0(\varepsilon
, k) \geq 1$ and a real $B = B(\varepsilon , k) > 0$ such that the
sequence ${([B. {n!}^{k + \varepsilon}])}_{n \geq n_0}$ is an
increasing sequence of prime numbers.
\end{coll}
\vskip 1mm \noindent {\bf Proof:} let $\epsilon$ be a fixed
positive real. From Conjecture \ref{conj1} (admitted with the
constant $k > 1$), there exists a positive real $c_k$ for which we
have:
\begin{equation}\label{eq7}
p_{n + 1} - p_n ~\leq~ c_k (\log{p_n})^k \hspace{2cm} (\forall n
\in \mathbb{N}) .
\end{equation}
We apply Theorem \ref{t1} for $I = ]1 , 2[$, $n_0 \geq 2$ an
integer (depending on $k$ and $\varepsilon$) which we pick large
enough that we have:
\begin{equation}\label{eq8}
c_k\left((k + \varepsilon) (n + 1) \log(n + 1) + \log{2}\right)^k
+ 1 ~\leq~ (n + 1)^{k + \epsilon} \hspace{1cm} (\forall n \geq
n_0) ,
\end{equation}
$f_n(x) = {n!}^{k + \varepsilon} x$ $(\forall n \geq n_0 , \forall
x \in I)$, $g(x) = c_k (\log{x})^k + 1$ if $x > 1$ and $g(x) = 1$
if $x \leq 1$ and ${(u_n)}_n$ the sequence of the prime numbers.
In this situation, we can easily check that the hypothesis of
Theorem \ref{t1} are all satisfied. We just note that the relation
(\ref{eq1}) follows from (\ref{eq8}), the relation (\ref{eq2})
follows from (\ref{eq7}) and the last hypothesis of Theorem
\ref{t1} concerning the sequence ${(u_n)}_n = {(p_n)}_n$ is a
consequence of Bertrand's postulate. Corollary \ref{coll2} follows
from this
application.\penalty-20\null\hfill$\blacksquare$\par\medbreak
\noindent

Apart from the context of the prime numbers, we have the following
\begin{coll}\label{coll3}
Let ${(u_n)}_{n \in \mathbb{N}}$ be a sequence of integers such
that: $\displaystyle 1 \leq \limsup_{n \rightarrow + \infty} (u_{n
+ 1} - u_n) < + \infty$. Then, we have:
\begin{enumerate}
\item[\rm{(1)}] For any positive real $\lambda$, there exists a
real $A > 1$ for which the sequence ${([\lambda A^n])}_{n \geq 1}$
is an increasing subsequence of ${(u_n)}_n$.
 \item[\rm{(2)}] For any real $A > \limsup_{n \rightarrow + \infty}
 (u_{n + 1} - u_n) + 1$, there exists a positive real $\lambda$ for which the sequence ${([\lambda
A^n])}_{n \geq 1}$ is an increasing subsequence of ${(u_n)}_n$.
\end{enumerate}
\end{coll}\vskip 1mm \noindent
{\bf Some open problems related to the preceding study:} We ask
(under or without Cramér's Conjecture) the following questions:
\begin{enumerate}
\item[(1)] Does there exists a real number $A > 1$ for which
$[A^n]$ is a prime number for any positive integer $n$? (This
corresponds to the case $\xi = 1$ which is excluded from Corollary
\ref{coll1}).
 \item[(2)] More generally than (1), does there exists a couple of real numbers $(\lambda , A)$, with
 $\lambda > 0$, $A > 1$, for which $[\lambda A^n]$ is a prime number for any positive integer $n$? (This is related
 to Corollary \ref{coll3}).
 \item[(3)] Does there exists a real number $B > 1$ for which $[B. {n!}^2]$ is a prime number
 for any sufficiently large non-negative integer $n$? (This is related to Corollary \ref{coll2}).
\end{enumerate}

\end{document}